\documentstyle[amssymb,draft,11pt]{amsart}
\title[Courant algebroids and SHLA's]{Courant algebroids and strongly homotopy 
Lie algebras}
 
\thanks{${}^*$The authors were supported in part by NSF grant DMS-9625122.}

\author[Dmitry Roytenberg and Alan Weinstein]{Dmitry Roytenberg\and Alan 
Weinstein
 ${}^*$
}

\address{Department of Mathematics, University of California, Berkeley, CA 
94720}
\email{roytenbe@@math.berkeley.edu}
\email{alanw@@math.berkeley.edu}
\date{\today}

\theoremstyle{plain}
\newtheorem{definition}{Definition}[section]
\newtheorem{theorem}[definition]{Theorem}
\newtheorem{lemma}[definition]{Lemma}
\newtheorem{example}[definition]{Example}
\newtheorem{prop}[definition]{Proposition}
\newtheorem{corollary}[definition]{Corollary}

\newcommand{\br}{[\cdot,\cdot ]}
\newcommand{\bform}[2]{\langle #1 , #2 \rangle}
\newcommand{\cinf}{C^{\infty}}
\newcommand{\lon }{\longrightarrow }
\newcommand{\half}{\frac{1}{2}}
\newcommand{\third}{\frac{1}{3}}
\newcommand{\fourth}{\frac{1}{4}}
\newcommand{\sixth}{\frac{1}{6}}
\newcommand{\smalcirc}{\mbox{\tiny{$\circ $}}}
\newcommand{\DD}{{\cal D}}
\newcommand{\coker}{\mbox{coker}\;}
\newcommand{\w}{\wedge}
\newcommand{\bw}{\bigwedge}

\begin{document}

\begin{abstract}
Courant algebroids are structures which include as examples the doubles of Lie 
bialgebras and the bundle $TM\oplus T^*M$ with the bracket introduced by T. 
Courant for the study of Dirac structures. Within the category of Courant 
algebroids one can construct the doubles of Lie bialgebroids, the infinitesimal 
objects for Poisson groupoids. We show that Courant algebroids can be considered 
as strongly homotopy Lie algebras.
\end{abstract}

\maketitle

\section{Introduction}
Dirac structures on manifolds were introduced and studied by T. Courant in his 
1990 paper \cite{Cou}. These structures provide a geometric setting for Dirac's 
theory of constrained mechanical systems \cite{Dir} in the same way as 
symplectic or Poisson structures do for unconstrained ones. A Dirac structure on 
a manifold $M$ is a subbundle $L\subset TM\oplus T^*M$ 
that is maximally isotropic with respect to the canonical symmetric bilinear 
form on $TM\oplus T^*M$, and which satisfies a certain integrability condition. 
To 
formulate the integrability condition, Courant introduced a skew-symmetric 
bracket operation on sections of $TM\oplus T^*M$; the condition is 
that the sections of $L$ be closed under this bracket. The Courant bracket 
is 
natural in the sense that it does not depend on any additional structure on $M$ 
for its definition, but it has anomalous properties. In particular, it does not 
satisfy the Leibniz rule with respect to multiplication 
by functions or the Jacobi identity. The ``defects'' in both cases are 
differentials of certain 
expressions depending on the bracket and the bilinear form; hence they disappear 
upon restriction to a Dirac subbundle because of the isotropy condition. 
Particular cases of Dirac subbundles are (graphs of) 2-forms and bivector 
fields on $M$, in which case the integrability condition turns out to coincide 
with the 
form being closed (resp. the bivector field being Poisson).

The nature of the Courant bracket itself remained unclear until several 
years later when it was observed by Liu, Weinstein and Xu \cite{LWX1} that 
$TM\oplus T^*M$ endowed with the Courant bracket plays the role of a ``double'' 
object, in the sense of Drinfeld (cf.\cite{Dr} or the book \cite{ChPr}), for 
a pair of {\em Lie algebroids} over $M$. Lie algebroids are structures on 
vector bundles that combine the features of both Lie algebras and vector fields. 
They, as well as the corresponding global objects, Lie groupoids, have recently 
found many applications in differential geometry \cite{Mac}, symplectic and 
Poisson geometry \cite{We1},\cite{We2},\cite{WeXu}, and also algebraic 
geometry and representation theory \cite{BeBe}. Many constructions in the 
category of Lie algebras carry over to Lie algebroids. Thus, in complete analogy 
with Drinfeld's Lie bialgebras, in the category of Lie algebroids there also 
exist ``bi-objects'', {\em Lie bialgebroids}, introduced by Mackenzie and Xu 
\cite{MacXu} as linearizations of Poisson groupoids. Besides their role in 
Poisson geometry and quantization \cite{KS4},\cite{We1},\cite{We2},\cite{WeXu}, 
Lie 
bialgebroids and Poisson groupoids have turned out to be the geometric 
structures behind the classical dynamical Yang-Baxter equation \cite{EtVar}.

It is well known \cite{ChPr},\cite{Dr} that every Lie bialgebra has a double 
which is a Lie algebra. This is not so for general Lie bialgebroids. 
Instead, Liu, Weinstein and Xu \cite{LWX1} show that the double of a Lie 
bialgebroid is a more complicated structure they call a {\em Courant algebroid}, 
$TM\oplus T^*M$ with the Courant bracket being a special case. On the other 
hand, when the base manifold is a point, all anomalies disappear and we get a 
Lie algebra with an invariant inner product.

The purpose of this work is to understand the anomalies of Courant algebroids. 
 From the very beginning, the observation that these anomalies were, in some 
sense, differentials of certain expressions suggested a homological/homotopical 
algebraic approach. This idea has turned out to be correct: we show that 
Courant algebroids can be viewed as 
{\em strongly homotopy Lie algebras} (also known as SHLA's or 
$L_\infty$-algebras). Homotopy Lie algebras play an important role in 
deformation 
theory as deformations of (differential, graded) Lie algebras 
\cite{Ko},\cite{ShlStash}, and they have also arisen in the latest 
developments in theoretical physics, in particular in string theory \cite{Zwi}; 
this is not a coincidence, since it has been shown \cite{Stash} that what 
physicists are dealing with is precisely deformation theory. For another 
appearance of homotopy Lie algebras in the context of constrained Hamiltonian 
mechanics see \cite{HenTeit}, or \cite{Kje} for a completely algebraic 
treatment.

Given a Courant algebroid, we realize it as a (finite) homological resolution of 
a Lie algebra. From general considerations it is natural to expect the 
resolution spaces of Lie algebras to carry homotopy Lie algebra structures 
(since a 
resolution of an object is an equivalent object in the derived category), and in 
\cite{BFLS} it is shown that a Lie algebra structure on a vector space can 
indeed be lifted to a SHLA structure on the total space of its homological 
resolution. However, the lifting depends on many choices, and it is not clear in 
advance how to make the ``best'' choice. In case of a Courant algebroid, we find 
explicit expressions for the structure maps which are natural and simple (in 
fact, most of them vanish), and prove that they satisfy the structure identities 
of a strongly homotopy Lie algebra.

The paper is organized as follows. In Section 2 we recall the basic notions 
related to Lie algebroids and bialgebroids; in Section 3 we introduce Courant 
algebroids and recall the construction of \cite{LWX1} of doubles of Lie 
bialgebroids; in Section 4 we introduce homotopy Lie algebras and prove the main 
theorem; Section 5 is devoted to the proofs of some technical lemmas; and 
Section 6 contains a few concluding remarks.

The authors would like to thank A.Nijenhuis, A. Schwarz, and J. Stasheff for 
sending us their papers, and for their encouraging comments.

\section{Lie algebroids and bialgebroids}

To begin, we recall the definition of Lie algebroid \cite{Mac}:
\begin{definition}
\label{def:liealg}
A {\em Lie algebroid} is a vector bundle $A\to M$ together with a Lie algebra 
bracket $\br$ on the space of sections $\Gamma(A)$ and a bundle map $\rho:A\to 
TM$, called the {\em anchor}, satisfying the following conditions:
\begin{enumerate}
\item For any $a_1, a_2\in\Gamma(A)$,
$\rho [a_1, a_2]=[\rho a_1,\rho a_2]$
\item For any $a_1, a_2\in\Gamma(A), f\in\cinf(M)$,
$[a_1, fa_2]=f[a_1, a_2]+(\rho(a_1)f)a_2$
\end{enumerate}
\end{definition}

In other words, the sections of the bundle act on smooth functions by 
derivations via the anchor in such a way that brackets act as commutators, 
and the behavior of the bracket with respect to multiplication by functions is 
governed by the Leibniz rule. Thus, Lie algebroids are a straightforward 
generalization of the tangent bundle. They are also the infinitesimal objects 
corresponding to Lie groupoids \cite{Mac}; when the base manifold is a point, a 
Lie groupoid reduces to a Lie group, while a Lie algebroid is just a Lie 
algebra.

A Lie algebroid structure on $A\to M$ gives rise to the 
following structures, dual to one another: first, the Lie bracket on $\Gamma(A)$ 
and the action of $\Gamma(A)$ on functions can be uniquely extended to a graded 
Lie algebra bracket of degree -1 on $\Gamma(\wedge^*A)$ which is a derivation of 
the exterior multiplication in each argument. This bracket is called the {\em 
Schouten bracket}, by analogy with the well-known bracket of multivector fields, 
and the resulting structure is a type of graded Poisson algebra called a {\em 
Gerstenhaber algebra}. 
Dually, one gets a differential $d_A$ on the graded commutative algebra 
$\Gamma(\wedge^*A^*)$, defined by the same formula as the usual de Rham 
differential and satisfying similar properties. The space $\Gamma(\wedge^*A^*)$ 
thereby acquires the structure of a differential graded commutative algebra.

Now suppose that we are given a pair $(A,A^*)$ 
of Lie algebroids over $M$ which are in duality as vector bundles. Then the Lie 
algebroid structure of $A$ induces a Schouten bracket on 
$\Gamma(\wedge^*A)$ and a differential $d_A$ on $\Gamma(\wedge^*A^*)$; on the 
other hand, from $A^*$ we get a Schouten bracket on $\Gamma(\wedge^*A^*)$ and 
a differential $d_{A^*}$ on $\Gamma(\wedge^*A)$.

\begin{definition}
\label{def:liebialg}
A pair $(A, A^*)$ of Lie algebroids in duality is a {\em Lie bialgebroid} if the 
induced differential $d_*$ is a derivation of the Schouten bracket on 
$\Gamma(\wedge^*A)$.
\end{definition}

It can be shown that this notion is self-dual (cf. Corollary \ref{cor} of the 
next section). Lie bialgebroids correspond to {\em differential Gerstenhaber 
algebras} \cite{KS1}. 

\begin{example}
Let $M$ be a Poisson manifold with Poisson tensor $\pi$ and the corresponding 
bundle map $\tilde{\pi}:T^*M\to TM$ given by $\langle\tilde{\pi}\alpha, 
\beta\rangle=\pi(\alpha,\beta)$. Let $A=TM$, the tangent bundle Lie algebroid,
$A^*=T^*M$  with anchor $\tilde{\pi}$ and the bracket of 1-forms given by
\begin{equation}
[\alpha,\beta]={\cal L}_{\tilde{\pi}\alpha}\beta-{\cal 
L}_{\tilde{\pi}\beta}\alpha-d(\pi(\alpha,\beta))
\end{equation}
Then $d$ is the usual deRham differential of forms, $d_*=[\pi,\cdot]$, and it is 
straightforward to verify that $(A, A^*)$ is a Lie bialgebroid.
\end{example}

Detailed discussion and more examples of Lie bialgebroids and Gerstenhaber 
algebras 
 from geometry and physics can be found in \cite{KS1},\cite{KS2} and \cite{KS4}. 

\section{Courant algebroids}
\begin{definition}
\label{def:jac}
Given a bilinear, skew-symmetric operation $\br$ on a vector space $V$, its {\em 
Jacobiator} $J$ is the trilinear operator on $V$:
\[
J(e_1,e_2,e_3)=[[e_1,e_2],e_3]+[[e_2,e_3],e_1]+[[e_3,e_1],e_2],
\]
$e_1,e_2,e_3\in V$.
\end{definition}
The Jacobiator is obviously skew-symmetric. Of course, in a Lie 
algebra $J\equiv 0$.
\begin{definition}
\label{def:quasi-algebroid}
A {\em Courant algebroid} is a vector bundle $E\lon M$
equipped with a nondegenerate  symmetric  bilinear form
 $\bform{\cdot}{\cdot}$ on the bundle,  a  skew-symmetric
bracket $\br$ on $\Gamma (E)$,
and a bundle map $\rho :E\lon TM$ such that the following
properties are  satisfied:
 \begin{enumerate}
\item For any $e_{1}, e_{2}, e_{3}\in \Gamma (E)$,
$J(e_1,e_2,e_3)={\cal D} T(e_{1}, e_{2}, e_{3});$
\item  for any $e_{1}, e_{2} \in \Gamma (E)$,
$\rho [e_{1}, e_{2}]=[\rho e_{1}, \rho  e_{2}];$
\item  for any $e_{1}, e_{2} \in \Gamma (E)$ and $f\in \cinf (M)$,
$[e_{1}, fe_{2}]=f[e_{1}, e_{2}]+(\rho (e_{1})f)e_{2}-
\langle e_{1}, e_{2}\rangle{\cal D} f ;$
\item $\rho \smalcirc {\cal D} =0$, i.e.,  for any $f, g\in \cinf(M)$,
$\langle{\cal D} f,  {\cal D}  g\rangle=0$;
\item for any $e, h_{1}, h_{2} \in \Gamma (E)$,
  $\rho (e) \langle h_{1}, h_{2}\rangle=\langle [e , h_{1}]+{\cal D} \langle e 
,h_{1}\rangle ,
h_{2}\rangle+\langle h_{1}, [e , h_{2}]+{\cal D}  \langle e 
,h_{2}\rangle\rangle$,
\end{enumerate}
where $T(e_{1}, e_{2}, e_{3})$ is the function on the base $M$
defined by:
\begin{equation}
\label{eq:T0}
 T(e_{1}, e_{2}, e_{3})=\third \langle[e_{1}, e_{2} ], e_{3}\rangle+c.p.,
\end{equation}
(``c.p.'' denotes the cyclic permutations of the $e_i$'s) and
${\cal D} :  \cinf(M)\lon \Gamma (E)$
is  the map   defined by ${\cal D} = \half \beta^{-1}\rho^{*} d_{0}$, 
 where $\beta $ is
the isomorphism between $E$
and $E^*$ given by the bilinear form and $d_0$ is the deRham differential. In 
other words,
\begin{equation}
\label{eq:D}
\langle{\cal D} f , e\rangle= \half  \rho (e) f .
\end{equation}
\end{definition}

In a Courant algebroid $E$, a {\em Dirac structure}, or {\em Dirac subbundle}, 
is a subbundle $L$ that is maximally isotropic under $\bform{\cdot}{\cdot}$ and 
whose sections 
are closed under $\br$. It is immediate from the definition that a Dirac 
subbundle is a Lie algebroid under the restrictions of the bracket and anchor.

Suppose now that both $A$ and $A^{*}$ are  Lie algebroids over the base
manifold $M$,  with anchors $a$ and $a_{*}$ respectively.
Let $E$ denote their   vector bundle direct sum:
$E=A\oplus A^{*}$.
On   $E$, there exist  two natural nondegenerate
bilinear forms, one symmetric and another antisymmetric:

\begin{equation}
\label{eq:pairing}
(X_{1}+\xi_{1} , X_{2}+\xi_{2})_{\pm}=\half (\langle \xi_{1},  X_{2}  \rangle
  \pm \langle \xi_{2} ,  X_{1}\rangle ).
\end{equation}

On $\Gamma (E)$, we  introduce a bracket  by

\begin{equation}
\label{eq:double}
[e_{1}, e_{2}]=([X_{1}, X_{2}]+L_{\xi_{1}}X_{2}-L_{\xi_{2}}X_{1}-d_{*}(e_{1},
 e_{2})_{-})
+ ([\xi_{1} , \xi_{2}]+L_{X_{1}}\xi_{2}-L_{X_{2}}\xi_{1} +d(e_{1}, e_{ 2})_{-}),
\end{equation}
where $e_{1}=X_{1}+\xi_{1}$ and $e_{2}=X_{2}+\xi_{2}$.\\\\

Finally, we let $\rho : E\lon TM$ be the bundle map defined by
$\rho =a +a_{*}$. That is,
\begin{equation}
\rho (X+\xi )=a(X)+a_{*} (\xi  ) , \ \ \forall X\in \Gamma (A) \mbox{ and }
\xi \in \Gamma (A^{*})
\end{equation}
It is easy to see that in this case the operator ${\cal D}$ as
defined by Equation (\ref{eq:D}) is given by
$${\cal D}=d_{*}+d, $$
where $d_{*}: \cinf(M)\lon \Gamma (A)$ and
$d: \cinf(M)\lon \Gamma ( A^{*})$ are the
usual differential  operators associated to Lie algebroids (cf. Sec. 2 and 
\cite{MacXu} for more details).

The following results, proved in \cite{LWX1}, show that the notion of 
Courant algebroid generalizes of the double construction 
to Lie bialgebroids:

\begin{theorem}
\label{thm:main1}
If  $(A, A^{*})$ is a Lie bialgebroid, then
$E=A \oplus A^{*}$ together with $(\br, \rho , (\cdot , \cdot)_{+})$
is a Courant algebroid.
\end{theorem}

\begin{theorem}
\label{thm:main2}
In a Courant algebroid $(E, \rho , \br,  \bform{\cdot}{\cdot})$,
suppose that  $L_{1}$ and $L_{2}$ are  Dirac subbundles  transversal
to each other, i.e., $E=L_{1}\oplus L_{2}$.
Then, $(L_{1}, L_{2} )$ is a Lie bialgebroid, where $L_{2}$  is
considered as the dual bundle of $L_{1}$  under the
pairing  $2\bform{\cdot}{\cdot}$.
\end{theorem}

An immediate consequence of the theorems above
is the following duality  property of  Lie bialgebroids,
which was first proved in \cite{MacXu}
and then by Kosmann-Schwarzbach \cite{KS1} using a simpler method.

\begin{corollary}
\label{cor}
If $(A, A^{*})$  is a Lie bialgebroid,  so is $(A^{*}, A)$.
\end{corollary}
\begin{example}
Given a manifold $M$, consider $TM$ with its standard Lie algebroid structure 
and $T^*M$ with zero anchor and bracket. Then $(TM, T^*M)$ is a Lie 
bialgebroid, and the double bracket (\ref{eq:double}) reduces to
\[
[X_1+ \xi_1, X_2 +\xi_2]=
[X_1,X_2]+ (L_{X_1}\xi_2-L_{X_2}\xi_1 +
d({\textstyle{\frac 12}}(\xi_1(X_2)-\xi_2(X_1))).
\]
This is the bracket originally introduced by Courant in \cite{Cou}. The anchor 
$\rho$ in this case is the projection to $TM$, and $\DD=d$, the deRham 
differential.
\end{example}
\begin{example}
When $M$ is a point, $(A,A^*)$ is a Lie bialgebra and $E$ is the usual Drinfeld 
double.
\end{example}

\section{Strongly homotopy Lie algebras and Courant algebroids}

Let $V$ be a graded vector space. Let $T(V)$ denote the tensor algebra of $V$ in 
the category of graded vector spaces, and let $\bw(V)$ denote its exterior 
algebra in the same category ($\bw(V)=T(V)/<v\otimes w+(-1)^{{\tilde v}{\tilde 
w}}w\otimes v>$, 
where ${\tilde v}$ denotes the degree of $v$). $T(V)$ (resp. $\bw(V)$) is 
not only an associative algebra, but also a coassociative coalgebra: if $V$ is 
of finite type, the comultiplication on $T(V)$ (resp. $\bw(V)$) is the 
adjoint of the multiplication on $T(V^*)$ (resp. $\bw(V^*)$), but in fact 
one does not need the dual space to define the comultiplication (see \cite{SHLA} 
for details).
\begin{definition}
\label{def:SHLA}
A {\em strongly homotopy Lie algebra} (SHLA, $L_\infty$-algebra) is a graded 
vector space $V$ together 
with a collection of linear maps
$l_k:\bw^k V\to
V$ of degree $k-2$, $k\geq 1$, satisfying the following relation for each 
$n\geq 1$ and for all homogeneous $x_1,\dots ,x_n\in V$:

\begin{equation}
\label{eqn:SHLA}
\sum_{i+j=n+1}(-1)^{i(j-1)}\sum_{\sigma}(-1)^{\sigma}\epsilon(\sigma)l_{j}(l_{i}
(x_{\sigma(1)},\dots ,x_{\sigma(i)}),x_{\sigma(i+1)},\dots ,x_{\sigma(n)})=0,
\end{equation}
where $\epsilon(\sigma)$ is the {\em Koszul sign} (arising from the fundamental 
convention of 
supermathematics that a minus sign is introduced whenever two consecutive 
odd elements are permuted), and $\sigma$ runs over all $(i,n-i)$-unshuffles 
(permutations satisfying $\sigma(1)<\dots <\sigma(i)$ and $\sigma(i+1)<\dots 
<\sigma(n)$) with $i\geq 1$.
\end{definition}

For $n=1$ this means simply that
$l_1$ is a differential on $V$; for $n=2$, $l_2$ is a 
superbracket on $V$ of 
which $l_1$ is a derivation (equivalently, $l_2:\bw^2(V)\to V$ is a 
chain map of complexes); $n=3$ gives the Jacobi identity for $l_2$ satisfied up 
to chain homotopy given by $l_3$, and higher $l_k$'s can be interpreted as 
higher homotopies. The algebraic theory of $L_\infty$-algebras is studied in 
\cite{HinSch} and \cite{SHLA} .

We shall write the equation (\ref{eqn:SHLA}) in the more succinct equivalent 
form:
\begin{equation}
\label{eqn:SHLA'}
\sum_{i+j=n+1}(-1)^{i(j-1)}l_jl_i=0,
\end{equation}
where we have extended each $l_i$ to all of $\bw(V)$ as a 
coderivation of the coalgebra structure on $\bw(V)$. This accounts for the 
permutations and signs in (\ref{eqn:SHLA}).

We are interested in $L_{\infty}$-algebras for the following reason: it is shown 
in \cite{BFLS} that, given a resolution $(X_*,d)$ of a vector space 
$H$ (graded or not), any Lie algebra structure on $H$ can be lifted to an 
$L_{\infty}$-algebra structure on the total resolution space $X$ with $l_1=d$.
The starting point of this construction is the observation that Lie brackets on 
$H$ correspond to bilinear skew-symmetric brackets on $X_0$ for which the 
boundaries form an ideal and the Jacobi identity is satisfied up to a boundary. 
This correspondence is in no way unique or canonical, as it 
requires a choice of a homotopy inverse to the quasi-isomorphism $(X_*,d)\to 
(H,0)$). But it is this latter bracket on $X_0$ that provides the starting point 
for constructing the SHLA structure on $X$, hence, if it is given, no choice is 
required at this stage, and we need never mention $H$. We shall presently see 
that with Courant algebroids we are in precisely this situation.

Let $E$ be a Courant algebroid over a manifold $M$. We know from the definition 
that the Courant bracket on $\Gamma(E)$ satisfies Jacobi up to a $\DD$-exact 
term. It turns out that, moreover, $Im(\DD)$ is an ideal in 
$\Gamma (E)$ with respect to the bracket. More precisely, the following identity 
holds: 
\begin{prop} 
For any $e\in\Gamma(E)$, $f\in C^{\infty}(M)$ one has
$$[e,\DD f]=\DD\bform{e}{\DD f}$$
\label{prop:main}
\end{prop}

\begin{pf*}{Proof}
Use axiom 5 in the definition of Courant algebroid with $e={\cal D}f$ and 
arbitrary $h_1$ and $h_2$, and then cyclically permute $e$, $h_1$ and $h_2$:

\begin{eqnarray*}
\rho(\DD f)\bform{h_1}{h_2}&=&\bform{[\DD f,h_1]+\DD\bform{\DD f}{h_1}}{h_2}+
                              \bform{h_1}{[\DD f,h_2]+\DD\bform{\DD f}{h_2}}\\
\rho(h_1)\bform{h_2}{\DD f}&=&\bform{[h_1,h_2]+\DD\bform{h_1}{h_2}}{\DD f}+
                              \bform{h_2}{[h_1,\DD f]+\DD\bform{h_1}{\DD f}}\\
\rho(h_2)\bform{\DD f}{h_1}&=&\bform{[h_2,\DD f]+\DD\bform{h_2}{\DD f}}{h_1}+
			      \bform{\DD f}{[h_2,h_1]+\DD\bform{h_2}{h_1}}.
\end{eqnarray*}
Now add the first two identities and subtract the third. Using Courant algebroid 
axioms 2, 4 and the definition of $\DD$, we get:
\[
\half\rho([h_1,h_2])f=\bform{\DD f}{2[h_1,h_2]}+\bform{h_1}{2[\DD f,h_2]}+
                      \bform{h_2}{2\DD\bform{\DD f}{h_1}}.
\]
Using the definition of $\DD$ again, we can rewrite this as:
\begin{eqnarray*}
0&=&\half\rho([h_1,h_2])f+\bform{h_1}{2[\DD f,h_2]}+\rho(h_2)\bform{h_1}{\DD 
f}=\\
 &=&\half\rho([h_1,h_2])f+\bform{h_1}{2[\DD              
                                           f,h_2]}+\half\rho(h_2)(\rho(h_1)f)=\\
 &=&\half\rho(h_1)(\rho(h_2)f)+\bform{h_1}{2[\DD f,h_2]} =\\
 &=&\bform{h_1}{\DD(\rho(h_2)f)+2[\DD f,h_2]}=\\
 &=&\bform{h_1}{2(\DD\bform{h_2}{\DD f}-[h_2,\DD f])}.
\end{eqnarray*}
The statement follows from the nondegeneracy of $\bform{\cdot}{\cdot}$.
\end{pf*}                   			                                                                 
\medskip
It will follow that we can extend the Courant bracket to an $L_\infty$-structure 
on 
the total space of the following resolution of $H=\coker\DD$:
\begin{equation}
\label{eqn:res}
\cdots\lon 0\lon X_2\stackrel{d_2}{\lon}X_1\stackrel{d_1}{\lon}X_0\lon H\lon 0,
\end{equation}
where $X_0=\Gamma(E)$, $X_1=\cinf(M)$, $X_2=\ker\DD$, $d_1=\DD$ and $d_2$ is the 
inclusion $\iota:\ker\DD\hookrightarrow\cinf(M)$. Remarkably, it turns out that, 
owing to the properties of Courant algebroids, the choices in 
the extension procedure can be made in a natural and simple way. 

Let us fix some notation: we will denote elements of $X_0$ by $e$, elements of 
$X_1$ by $f$ or $g$, and elements of $X_2$ by $c$.

\begin{theorem} 
\label{thm:main}
A Courant algebroid structure on a vector bundle $E\lon M$ gives rise naturally 
to a SHLA structure on the total space $X$ of (\ref{eqn:res}) with $l_1=d$ and 
the higher structure maps given by the following explicit formulas: 
\[
\begin{array}{lccl}
l_2(e_1\w e_2)&=&[e_1,e_2]&\mbox{in degree 0}\\
l_2(e\w f)    &=&\bform{e}{\DD f}&\mbox{in degree 1}\\
l_2	      &=&0&\mbox{in degree $> 1$}\\
l_3(e_1\w e_2\w e_3)&=&-T(e_1,e_2,e_3)&\mbox{in degree 0}\\
l_3                 &=&0&\mbox{in degree $> 0$}\\
l_n  &=&0&\mbox{for $n>3$}
\end{array}
\]
\end{theorem}

\begin{pf*}{Proof}
Starting with the Courant bracket on $X_0$, we shall, following \cite{BFLS}, 
extend it to an $l_2$ on all of $X$ satisfying (\ref{eqn:SHLA'}) for $n=2$. The 
extension will proceed, essentially, by induction on the degree of the argument: 
for each degree $l_2$ will be a primitive of a certain cycle depending on the 
values of $l_2$ on elements of lower degree. Higher $l_k$'s will be introduced 
and extended in a similar fashion, as primitives of cycles (using the acyclicity 
of (\ref{eqn:res})). The main work will consist in calculating these 
cycles, in particular, showing that most of them vanish; these computations are 
mostly relegated to the technical lemmas of the next section. 

\noindent{\em Step 1: $n=2$}. In degree 0, we are given $l_2(e_1\w e_2)=
[e_1,e_2]$. Consider now an element $e\w f$ of degree 1. Then $l_2l_1(e\w f)\in 
X_0$ 
is defined and is, in fact, a boundary by Proposition \ref{prop:main}:
\[
l_2l_1(e\w f)=l_2(l_1e\w f+e\w l_1f)=[e,\DD f]=\DD\bform{e}{\DD f},
\]
so we set $l_2(e\w f)=\bform{e}{\DD f}$ so that the SHLA identity 
(\ref{eqn:SHLA'}) for $n=2$,
\begin{equation}
\label{eqn:SHLA2}
l_1l_2-l_2l_1=0,
\end{equation} 
holds in degree 1.

Now, $\bw^2(X)_2$ is spanned by elements of the form $f\w g$ or $c\w e$. As 
above, $l_2l_1$ is defined on elements of degree 2, and is, in fact, a cycle 
(cf. \cite{BFLS}). We have
\[
l_2l_1(f\w g)=l_2(l_1f\w g-f\w l_1g)=l_2(\DD f\w g-f\w\DD g)=\bform{\DD f}{\DD 
g}+\bform{\DD g}{\DD f}=0
\]
by Courant algebroid axiom 4, whereas
\[
l_2l_1(c\w e)=l_2(l_1c\w e+c\w l_1e)=l_2(\iota c\w e)=-\bform{e}{\DD\iota c}=0,
\]
so we set $l_2(f\w g)=l_2(c\w e)=0$. Now observe that, since $l_2=0$ in degree 
2, we can define $l_2$ to be 
zero on elements of degree higher than 2 as well and still have 
(\ref{eqn:SHLA2}). We have thus defined an $l_2$ that satisfies 
(\ref{eqn:SHLA2}) by construction.

\noindent{\em Step 2: $n=3$}. In degree 0, by Courant algebroid axiom 
1 we have
\[
l_2l_2(e_1\w e_2\w e_3)=J(e_1,e_2,e_3)=\DD T(e_1,e_2,e_3),
\]
where $J$ is the Jacobiator. So we set 
$l_3(e_1\w e_2\w e_3)=-T(e_1,e_2,e_3)$, so that the homotopy Jacobi identity 
identity (\ref{eqn:SHLA'}) for $n=3$,
\begin{equation}
\label{eqn:SHLA3}
l_1l_3+l_2l_2+l_3l_1=0,
\end{equation}
holds on $\bw^3(X)_0$ (as $l_1(X_0)=0$). 

Consider now an element $e_1\w e_2\w f\in\bw^3(X)_1$ . The expression 
$(l_2l_2+l_3l_1)(e_1\w e_2\w f)$ is defined and is a cycle in $X_1$ (cf. 
\cite{BFLS}), hence we can define $l_3(e_1\w e_2\w f)$ to be some primitive of 
this cycle, so that (\ref{eqn:SHLA3}) holds. But in our particular situation we 
in fact have (see the next section for a proof): 
\begin{lemma}
\label{lemma:t1}
$(l_2l_2+l_3l_1)(e_1\w e_2\w f)=0$ $\forall e_1,e_2,f$.
\end{lemma}
Therefore, we can define $l_3(e_1\w e_2\w f)=0$. Now observe that on elements of 
degree $>1$ $l_3$ has to be 0 because deg$(l_3)=1$, whereas $X_k=0$ for $k>2$. 
We now have $l_3$ defined on all of $\bw^3(X)$ and satisfying (\ref{eqn:SHLA3}) 
by construction.

\noindent{\em Step 3: $n=4$ and higher}. Proceeding in a similar 
fashion, we look at the expression\\ $(l_3l_2-l_2l_3)(e_1\w e_2\w e_3\w e_4)$ 
(always a cycle in $X_1$) and define $l_4(e_1\w e_2\w e_3\w e_4)$ to be its 
primitive in $X_2$, so as to satisfy (\ref{eqn:SHLA'}). However, it turns out 
that (see the next section for a proof)
\begin{lemma}
\label{lemma:t2}
$(l_3l_2-l_2l_3)(e_1\w e_2\w e_3\w e_4)=0$ $\forall e_1,e_2,e_3,e_4$.
\end{lemma}
Hence we can set $l_4(e_1\w e_2\w e_3\w e_4)=0$ and observe that $l_4$ has to 
vanish on elements of degree $>0$ as deg$(l_4)=2$, while $X_k=0$ for $k>2$. By 
similar degree counting, all $l_n$, $n>4$, have to vanish identically. This 
finishes the proof 
modulo Lemmas \ref{lemma:t1} and \ref{lemma:t2}.
\end{pf*}

\section{Proofs of technical lemmas}
Let $(E,\bform{}{},\br,\rho)$ be a 
Courant algebroid over $M$. Given $e\in\Gamma(E)$, $f\in\cinf(M)$, we will 
denote $\rho(e)f$ simply by $ef$, for short. Let us first prove two auxiliary 
lemmas. 
\begin{lemma}
\label{lemma:a1} The identity
\[
T(e_1,e_2,\DD f)=\fourth[e_1,e_2]f
\]
holds in any Courant algebroid.
\end{lemma}
\begin{pf*}{Proof}
Using Courant algebroid axiom 2 and Proposition \ref{prop:main}, we have
\begin{eqnarray*}
T(e_1,e_2,\DD f)&=&\third(\bform{[e_1,e_2]}{\DD f}+\bform{[\DD f,e_1]}{e_2}+
\bform{[e_2,\DD f]}{e_1})=\\
&=&\third(\bform{[e_1,e_2]}{\DD f}-\bform{\DD\bform{e_1}{\DD f}}{e_2}+
\bform{\DD\bform{e_2}{\DD f}}{e_1})=\\
&=&\third(\half[e_1,e_2]f-\fourth e_2(e_1f)+\fourth e_1(e_2f))=\\
&=&\third(\half[e_1,e_2]f+\fourth[e_1,e_2]f)=\fourth[e_1,e_2]f.
\end{eqnarray*}
\end{pf*}

\begin{lemma}
\label{lemma:a2}
Given $e_1, e_2, e_3, e_4\in\Gamma(E)$, let
\begin{eqnarray*}
{\bf J}&=&\bform{J(e_1,e_2,e_3)}{e_4}-\bform{J(e_1,e_2,e_4)}{e_3}+
\bform{J(e_1,e_3,e_4)}{e_2}-\bform{J(e_2,e_3,e_4)}{e_1}\\
{\bf K}&=&\bform{[e_1,e_2]}{[e_3,e_4]}-\bform{[e_1,e_3]}{[e_2,e_4]}+
\bform{[e_1,e_4]}{[e_2,e_3]},
\end{eqnarray*}
where $J$ is the Jacobiator (cf. Def \ref{def:jac}). Then ${\bf K}+2{\bf J}=0$.
\end{lemma}
\begin{pf*}{Proof}
Using Courant algebroid axioms 1 and 5, we can rewrite ${\bf J}$ as follows:
\begin{eqnarray*}
\bform{J(e_1,e_2,e_3)}{e_4}&=&\bform{\DD T(e_1,e_2,e_3)}{e_4}=\half 
e_4T(e_1,e_2,e_3)=\sixth e_4(\bform{[e_1,e_2]}{e_3}+c.p.)=\\
&=&\sixth(\bform{[e_4,[e_1,e_2]]+\DD\bform{e_4}{[e_1,e_2]}}{e_3}+
\bform{[e_1,e_2]}{[e_4,e_3]+\DD\bform{e_4}{e_3}})+c.p.
\end{eqnarray*}
Expressing the other summands of ${\bf J}$ in this form and collecting like 
terms in the parentheses, we find that the terms of the form 
$\bform{[e_i,e_j]}{\DD\bform{e_k}{e_l}}$ cancel out, terms of the form 
$\bform{[e_i,e_j]}{[e_k,e_l]}$ add up to $-4{\bf K}$, those of the form 
$\bform{[e_i,[e_j,e_k]]}{e_l}$ add up to ${\bf J}$, and finally, terms of the 
form $\bform{\DD\bform{e_i}{[e_j,e_k]}}{e_l}$ add up to $-3{\bf J}$ after we 
use Courant algebroid axiom 1. Thus,
\[
{\bf J}=\sixth({\bf J}-3{\bf J}-4{\bf K}),
\]
and the statement of the lemma follows immediately.
\end{pf*}

\begin{pf*}{Proof of Lemma \ref{lemma:t1}} In the notation of the previous 
section, we have, using Lemma \ref{lemma:a1} and Courant algebroid axiom 2:
\begin{eqnarray*}
&&(l_2l_2+l_3l_1)(e_1\w e_2\w f)=\\
&=&l_2(l_2(e_1\w e_2)\w f+l_2(e_2\w f)\w e_1+l_2(f\w e_1)\w e_2)+\\
&+&l_3(l_1e_1\w e_2\w f+e_1\w l_1e_2\w f+e_1\w e_2\w l_1f)=\\
&=&l_2([e_1,e_2]\w f+\bform{e_2}{\DD f}\w e_1-\bform{\DD f}{e_1}\w 
e_2)+l_3(e_1\w e_2\w \DD f)=\\
&=&\bform{[e_1,e_2]}{\DD f}-\bform{e_1}{\DD\bform{e_2}{\DD f}}+
\bform{e_2}{\DD\bform{e_1}{\DD f}}-T(e_1,e_2,\DD f)=\\
&=&\half[e_1,e_2]f-\fourth
e_1(e_2f)+\fourth e_2(e_1f)-\fourth[e_1,e_2]f=0
\end{eqnarray*}
\end{pf*}

\begin{pf*}{Proof of Lemma \ref{lemma:t2}}
In the notation of the previous section we have
\begin{eqnarray*}
l_2l_3(e_1\w e_2\w e_3\w e_4)&=&
l_2(l_3(e_1\w e_2\w e_3)\w e_4\pm(3,1)-unshuffles)=\\
&=&-l_2(T(e_1,e_2,e_3)\w e_4\pm(3,1)-unshuffles)=\\
&=&\bform{\DD T(e_1,e_2,e_3)}{e_4}\pm(3,1)-unshuffles=\\
&=&\bform{J(e_1,e_2,e_3)}{e_4}\pm(3,1)-unshuffles={\bf J}.
\end{eqnarray*}
On the other hand,
\begin{eqnarray*}
l_3l_2(e_1\w e_2\w e_3\w e_4)&=&l_3(l_2(e_1\w e_2)\w e_3\w 
e_4)\pm(2,2)-unshuffles=\\
&=&-T([e_1,e_2],e_3,e_4)\mp(2,2)-unshuffles=\\
&=&-\third(\bform{[e_1,e_2],e_3]}{e_4}+\bform{[e_3,e_4]}{[e_1,e_2]}+
\bform{[e_4,[e_1,e_2]]}{e_3})\pm\cdots=\\
&=&-\third({\bf J}-2{\bf K}),
\end{eqnarray*}
after collecting like terms. An application of Lemma \ref{lemma:a2} 
immediately yields $l_2l_3=l_3l_2$.
\end{pf*}

\section{Concluding remarks}
$L_\infty$-algebras occur in physics in the framework of the Batalin-Vilkovisky 
procedure for quantizing gauge theories. On the other hand, the Courant bracket 
seems to provide a geometric framework for constrained Hamiltonian systems. It 
is known \cite{HenTeit} that gauge Lagrangians lead to constrained theories in 
the Hamiltonian formalism. This suggests that homotopy Lie algebras arising in 
the Batalin-Vilkovisky formalism and those in the Courant formalism might be 
somehow related. Our current investigations are aimed in this direction.


\providecommand{\bysame}{\leavevmode\hbox to3em{\hrulefill}\thinspace}

\end{document}